\documentclass[sn-mathphys]{sn-jnl}

\graphicspath{ {./Figures/} }
\usepackage{subfig} 
\usepackage{booktabs}
\usepackage{amsmath}
\usepackage{xcolor}
\usepackage{adjustbox}
\usepackage{multicol}
\usepackage{enumitem}
\usepackage{xcolor}
\usepackage{float}
\definecolor{MyBlue}{rgb}{0.25,0.5,0.75}
\colorlet{NextBlue}{MyBlue!20}
\colorlet{SecondBlue}{MyBlue!40}
\jyear{2021}

\begin{document}

\title{Experimental and numerical investigations on heat transfer in fused filament fabrication 3D-printed specimens}

\author*[1,2]{\fnm{Nathalie} \sur{Ramos}}\email{nathalie.ramos@ntnu.no}
\author[2]{\fnm{Christoph} \sur{Mittermeier}}
\author[2,1]{\fnm{Josef} \sur{Kiendl}}

\affil[1]{\orgdiv{Department of Marine Technology}, \orgname{Norwegian University of Science and Technology}, \orgaddress{\street{Otto Nielsens Veg 10}, \city{Trondheim}, \postcode{7491}, \country{Norway}}}

\affil[2]{\orgdiv{Institute of Engineering Mechanics \& Structural Mechanics}, \orgname{Bundeswehr University Munich}, \orgaddress{\street{Werner-Heisenberg-Weg 39}, \city{Neubiberg}, \postcode{85577}, \country{Germany}}}

\abstract{A good understanding of the heat transfer in fused filament fabrication is crucial for an accurate stress prediction and subsequently for repetitive, high quality printing. This work focuses on two challenges that have been presented when it comes to the accuracy and efficiency in simulating the heat transfer in the fused filament fabrication process. With the prospect of choosing correct thermal boundary conditions expressing the natural convection between printed material and its environment, values for the convective heat transfer coefficient and ambient temperature were calibrated through numerical data fitting of experimental thermal measurements. Furthermore, modeling simplifications were proposed for an efficient numerical discretization of infill structures. Samples were printed with varying infill characteristics, such as varying air void size, infill densities and infill patterns. Thermal measurements were performed to investigate the role of these parameters on the heat transfer and based on these observations, possible modeling simplifications were studied in the numerical simulations.}


\maketitle

\section{Introduction}
\label{sec:intro}

\begin{figure*}
\centering
\subfloat[Programmed bending]{\includegraphics[width=0.3\textwidth]{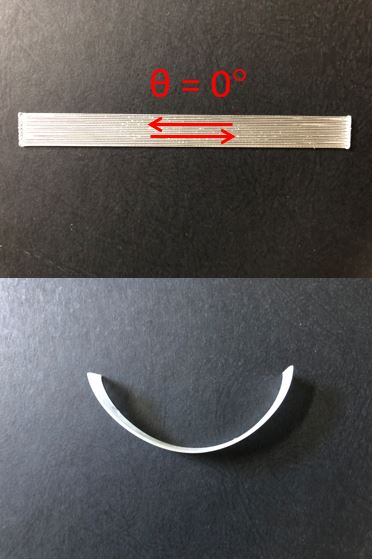}}\hspace{1em}
\subfloat[Programmed twisting]{\includegraphics[width=0.3\textwidth]{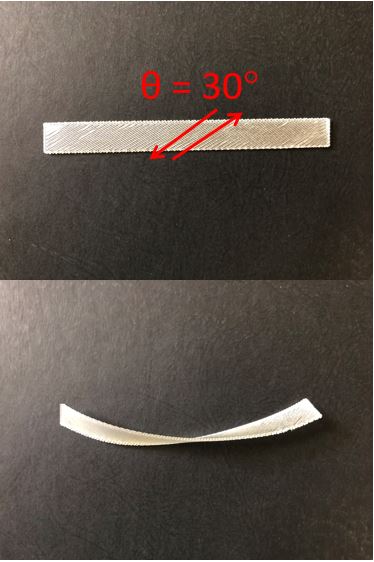}}
\caption{Shape change after heating FFF printed specimens above T$_g$} 
\label{fig:SMP}
\end{figure*}

Fused filament fabrication (FFF), or fused deposition modeling, is one of the most widely applied additive manufacturing methods. It is a three-dimensional (3D) printing method in which a thermoplastic material is extruded through a nozzle to construct a layer-by-layer structure \cite{Cattenone2019}. Once lauded for its potential in prototype manufacturing, it has now found its way into a plethora of fields and applications, such as biomedical, electrical and aerospace engineering \cite{Barry2012}. Advantages of FFF printing are the relatively low costs, wide applicability and availability, large variety of suitable materials and the ease of use \cite{Yang2016,Tymrak2014}. Some of these advantages have caused a surge of the use of FFF in four-dimensional (4D) printing, where the fourth dimension represents the change of shape over time of the smart materials \cite{Tibbits2014,Farid2021}. By combining FFF printing, and the use of shape memory polymers (SMP), structures can be printed that can maintain a temporary shape and that can return to their original shape after being exposed to an external stimulus, such as heat \cite{Li2017}. The thermal sensitive nature of printed SMP parts enables potential applications such as fasteners in active assembly/disassembly, smart actuators and deployable structures for aerospace applications  \cite{Yang2016,Ge2013,Ge2014}. \\ 
\indent However, deposition of the filament at high temperatures, followed by quick cooling to enforce solidification results in significant thermal gradients and subsequently, residual stresses \cite{Cattenone2019}. The presence of residual stresses can lead to warping during printing which can result in a failed print \cite{Fitzharris2018}. This is a common problem for Acrylonitrile Butadiene Styrene (ABS). Residual stresses can also cause distortion of parts and a loss of strength \cite{Kuznetsov2018}. In case of printed parts with strict geometrical tolerances or structural requirements, this effect can lead to a loss of functionality.\\
\indent During FFF printing the material is extruded at a temperature which is typically higher than the glass transition temperature. It is either deposited on the building platform or cooled existing layers, causing the material to be (partially) stretched as it bonds, to cool down and to finally solidify \cite{Hu2017}. The induced pre-strain will be released as soon as the material is heated above the glass transition temperature and a 'new' permanent shape emerges. An example of such programming can be seen in figure \ref{fig:SMP} where two rectangular specimens are printed flat, but with a different orientation of the printed filament. After reheating the samples above the glass transition temperature, they either bend upwards, or twist upwards depending on the printing orientation. Such approach has been used in 4D printing to design structures with self-folding, self-bending, self-twisting and shape-shifting mechanisms \cite{Hu2017,Bodaghi2017,Manen2017,Manen2018}.\\ 
\indent Whether the objective is to repetitively print parts of high quality, or to print SMP parts with a 4D effect, understanding the heat transfer is crucial for an accurate stress and bond strength prediction \cite{Yin2018}. Several efforts have been made to predict thermal gradients and the development of the residual stresses in printed parts by performing thermo-mechanical finite element simulations. Such simulations often make use of sequential element activation to represent the deposition sequence that takes place during FFF printing. Zhang and Chou \cite{Zhang2008} used such a model to perform a parametric study to predict part distortions in ABS. The analysis was thermo-mechanically coupled and the influence of various process parameters on the residual stresses was studied. Zhou et al.\cite{Zhou2016} presented a numerical model in which the heat transfer in ABS material was analyzed during the FDM process. In their model temperature-dependent material properties were used for the specific heat capacity and thermal conductivity. Cattenone et al. \cite{Cattenone2019} performed thermo-mechanical simulations in which the FFF process was simulated to predict part distortions in ABS. Special attention was paid to the constitutive modeling of the polymer, and the influence of numerical considerations on the predicted residual stresses, such as time step size and meshing strategies. Yin et al. \cite{Yin2018} performed similar heat transfer simulations but with a different objective. They used such models to predict the inter-facial bonding strength between printed filaments.\\ 
\indent This work focuses on two challenges that have been presented when it comes to the accuracy and efficiency in simulating the heat transfer in the FFF process. The first challenge that arises is the correct choice of thermal boundary conditions, particularly the convective heat transfer coefficient. The exact magnitude of this parameter is not always explicitly stated, its empirical determination can be quite cumbersome and often the focus is on forced convection.  Pereira et al. \cite{Pereira2019} focused on forced convection in their investigation of the effect of surface roughness on the convective heat transfer on the surfaces of FFF printed ABS cylindrical specimens. Zhou et al. \cite{Zhou2017} also calculated a value of the heat transfer coefficient based on forced convection over a rectangular body model. Costa et al. \cite{Costa2015} did focus on natural convection between printed filament and the ambient air, but the determined values for the heat transfer coefficient still varied between a fairly large range (5-60 W/m$^2\cdot$K) and this range was not experimentally validated. Lepoivre et al. \cite{Lepoivre2020} determined the heat transfer coefficient by an empirical correlation for an external free convection flow. The magnitude of the parameters used in said correlation were not listed or further elaborated.\\
\indent In this work the heat transfer coefficient is determined experimentally. Experimental thermal measurements are numerically simulated and a value for the heat transfer coefficient that describes natural convection is determined through data fitting.\\
\indent The second challenge is related to the computational effort which is required for an accurate heat transfer simulation. In previous work, the printed material is often modeled as a continuum in which the finite element size is determined by the filament cross-sectional dimensions. This is a computationally expensive exercise as the cross-section of a filament is much smaller than the global dimensions of printed geometries. As a result, the question arises how the modeling can be simplified to speed up simulations. The first question that must be answered is whether the assumption of discretizing the material as a continuum without accounting for the inherent air voids is an accurate one. If this is the case, methods to simplify the discretization of the characteristic mesostructure can be investigated and the computational effort required in heat transfer simulations can be reduced. Thus, in this work samples are printed with varying infill characteristics, such as varying air void size, infill densities and infill patterns. Thermal measurements are performed to investigate the role of these parameters on the heat transfer. Based on the observations made in these experiments, possible modeling simplifications are studied in the numerical simulations.\\
\indent The remainder of this paper is set up as follows. The materials and methods used in the thermal measurements and numerical simulations on samples of varying infill geometries are presented in section \ref{sec:mat_meth}. The results of the experimental thermal measurements are subsequently presented in section \ref{sec:exp_res}. Section \ref{sec:num_sim} is dedicated to the numerical simulations. A value for the heat transfer coefficient is determined by numerically fitting the data obtained in the experiments and modeling simplifications of complex infill geometries are investigated. Finally the conclusions are presented in section \ref{sec:conclu} and ideas for future research are proposed.

\section{Materials and Methods}
\label{sec:mat_meth}

\subsection{Experimental set-up}
\label{subsec:exp_setup}
\begin{figure*}
\centering
\subfloat[Set-up thermal recording]{\includegraphics[width=0.45\textwidth]{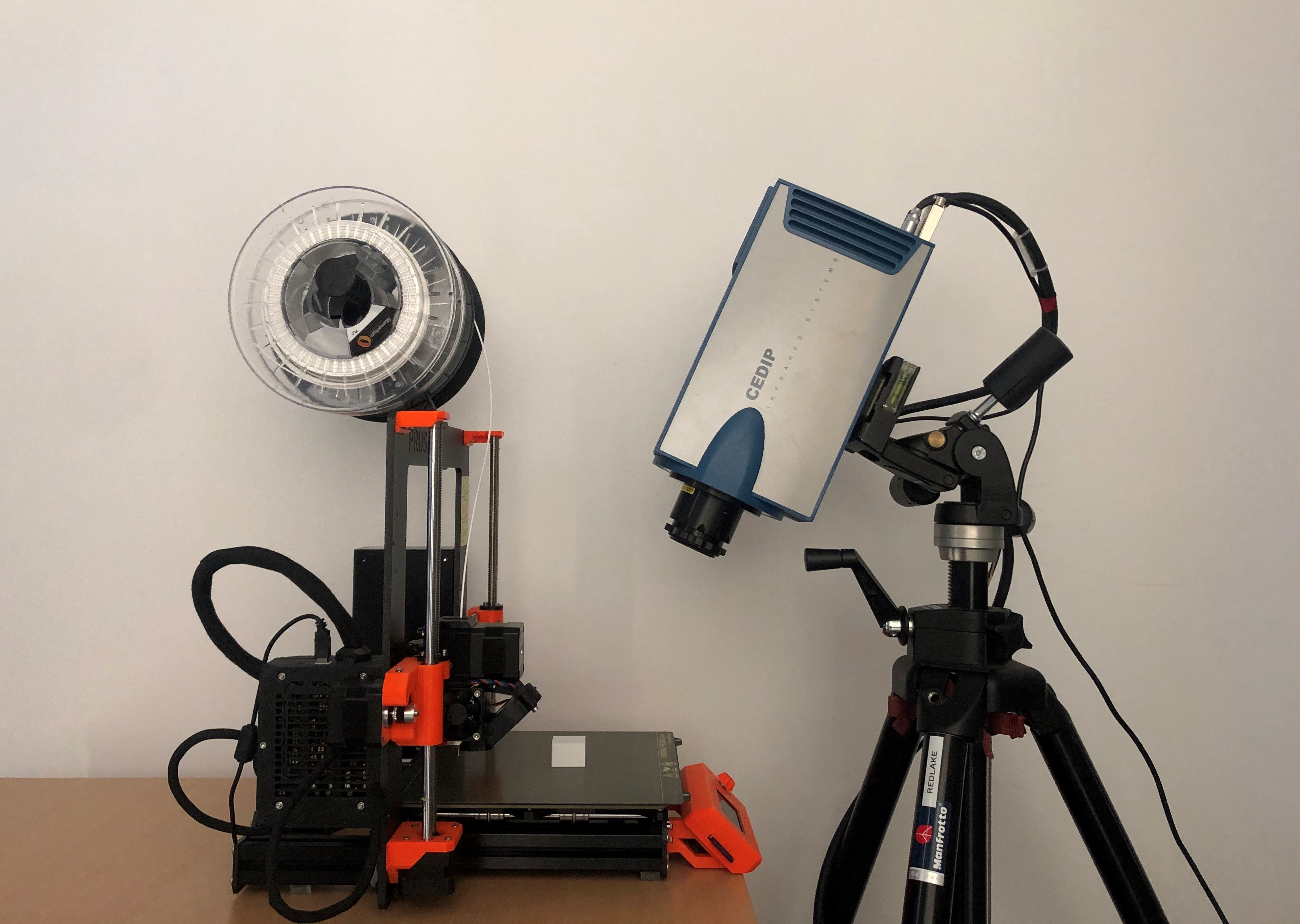}}\hfill
\subfloat[Data points thermal measurements]{\includegraphics[width=0.45\textwidth]{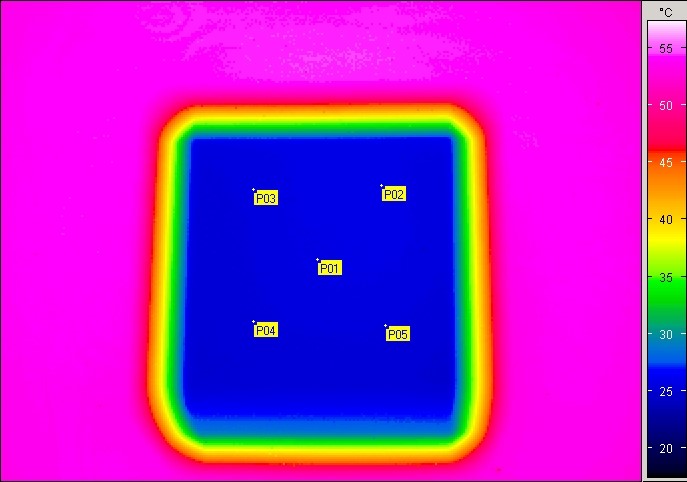}}
\caption{Experimental set-up}
\label{fig:setup}
\end{figure*}

All specimens used in this work were printed with a Prusa i3 MK3 printer. The material used was Poly Lactid Acid (PLA) and the thermal properties as provided by the manufacturer Fillamentum are listed in table \ref{tab:mat_prop}. A JADE CW infrared camera from Cedip Infrared Systems was used for the thermal measurements. The experimental set-up used in this paper is shown in figure \ref{fig:setup}. The measuring procedure was as follows. First the printing bed was heated up to a nominal temperature of 60$^{\circ}$C. A printed specimen was then placed on the heated printing bed for a sufficiently long time to allow for time to reach a steady-state. The heating of the specimens was recorded with the thermal camera and the temperature profiles were recorded with a sampling frequency of 1 Hz. The temperature evolution in time was recorded for five points on the top surface of the specimens (figure \ref{fig:setup}). All measurements were performed three times and they were spaced at least 4 hours apart to ensure full cooling prior to reheating. Even though the temperature of the printing bed was set to 60$^{\circ}$C, the temperature which was measured on the surface of the bed with the thermal camera was equal to 56$^{\circ}$C. The temperature of the ambient air in the room was 25$^{\circ}$C.

\begin{table} 
\centering
\begin{tabular}{lc}
\toprule 
\textbf{Property} & \textbf{Value}\\
\midrule
Density $\rho$ & 1240 kg/m$^{3}$ \\
Specific heat capacity $c_p$ & 1800 J/kg$\cdot$K \\
Conductivity $K_{0}$ & 0.13 W/m$\cdot$K\\
\bottomrule 
\end{tabular}
\caption{Thermal properties PLA}
\label{tab:mat_prop}
\end{table}

\subsection{Geometries and infill structure}
\label{subsec:geom_infill}

To facilitate the study of the two objectives in this paper different specimens printed with varying global and infill geometries were used during the thermal measurements. A reference specimen was printed which was used in the fitting process of the convective thermal boundary condition. All of the other printed specimens had a different infill geometry (air content, infill pattern) to investigate the influence of different infills on heat flow through the specimen. There are two different ways to vary the air content in the FFF printed specimens:

\begin{enumerate}
\item By variation of the extrusion factor
\item By variation of the infill density
\end{enumerate}

\begin{figure*}
\centering
\subfloat[Default infill]{\includegraphics[width=0.3\textwidth]{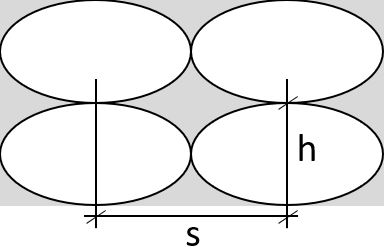}}\hfill
\subfloat[Varying the extrusion factor]{\includegraphics[width=0.3\textwidth]{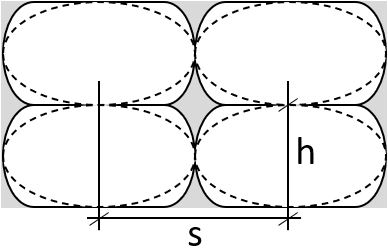}}\hfill
\subfloat[Varying the infill density]{\includegraphics[width=0.34\textwidth]{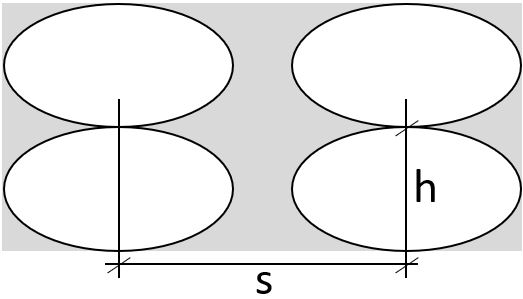}}
\caption{Two different ways of varying the air content}
\label{fig:air_content}
\end{figure*}

\begin{figure*}
\centering
\subfloat[CT-scan of the cross-section of a printed specimen]{\includegraphics[width=0.45\textwidth]{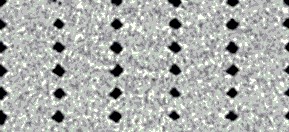}}\hspace{1em}
\subfloat[Schematization filaments \& air voids]{\includegraphics[width=0.35\textwidth]{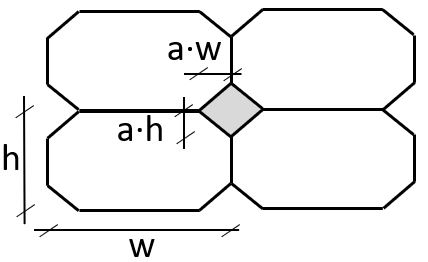}}\\
\subfloat[FE mesh of filaments and air voids]{\includegraphics[width=0.45\textwidth]{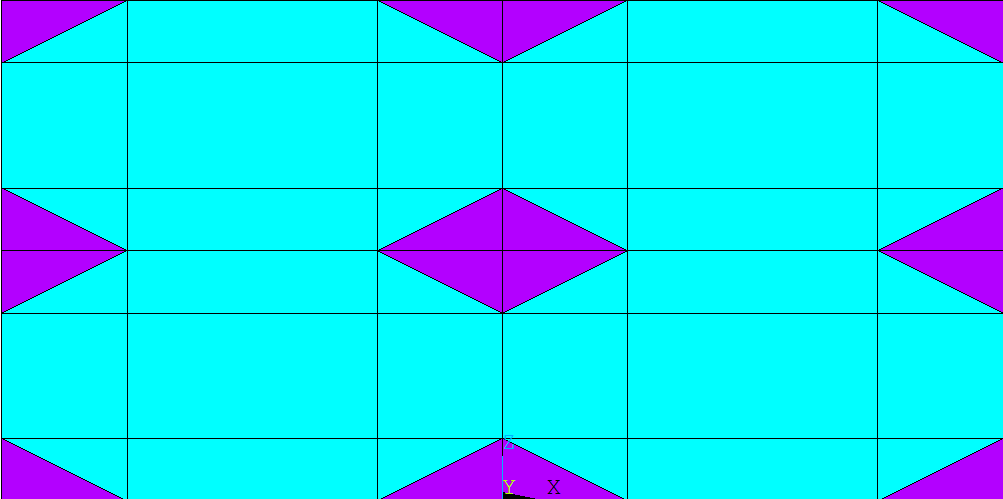}}
\caption{Discretization of air voids in the printed samples}
\label{fig:discretization}
\end{figure*}

\noindent Both approaches are shown schematically in figure \ref{fig:air_content}. Varying the extrusion factor in a densely packed setting of filaments influences the size of the air voids between those filaments. The infill density controls the gap size between the printed filaments.\\
\indent The Prusa slicer recommends a default value for the extrusion factor. This will be referred to as the default extrusion factor of 1.0 (EF=1.0). One sample was printed with the default extrusion factor, and the second sample was printed with an extrusion factor which was 10\% higher than the default value. The dimensions of these samples were 30x30x2 mm and both specimens were printed with an infill density of 100\%. These specimens were heated for five minutes. For the numerical discretization of the geometries with varying extrusion factor, a model as shown in figure \ref{fig:discretization} was used. The diamond shaped air voids were assumed based on a CT-scan of one of the printed samples (figure \ref{fig:discretization}). The air void size is governed by parameter $a$ and it was used to express the area of the air voids $A_{air}$ as a fraction of the total area $A_{tot}$: 

\begin{equation}
\label{eq:a}
\begin{aligned}
A_{tot}&=w\cdot{h}\\[2pt]
A_{air}&=\frac{1}{2}\cdot{aw}\cdot{ah}\cdot{4}=2\cdot{a}^{2}wh\\[2pt]
A_{f}&=A_{tot}-A_{air}=(1-2a^2)\cdot{wh}\\[2pt]
v_{fr;a}&=\frac{A_{air}}{A_{tot}}=2a^2
\end{aligned}
\end{equation}

\noindent where $w$ and $h$ are the width and the height of a printed filament respectively, as can be seen in figure \ref{fig:discretization}. The width of the filament equaled 0.45 mm and the layer height equaled 0.2 mm for all samples used in this work. By calculating the real volume fractions of air and filament in the printed samples, eq. \ref{eq:a} was solved for $a$. The real volume fractions were calculated as follows (eq. \ref{eq:volume_frac}):
\begin{enumerate}
\item The mass ($m$) and total volume ($V_{tot}$) of each printed sample were measured. The total volume was calculated by measuring the length, width and thickness of the printed specimens with a vernier caliper.
\item From the mass, the (real) extruded volume $V_{extr}$ was calculated. The density of PLA was used (table \ref{tab:mat_prop}). 
\item The volume fractions of the filament $v_{fr;f}$ and air $v_{fr;a}$ were then calculated.
\end{enumerate}

\begin{equation}
\label{eq:volume_frac}
\begin{aligned}
V_{extr}&=\frac{m}{\rho}\\[2pt]
v_{fr;f}&=\frac{V_{extr}}{V_{tot}}\\[2pt]
v_{fr;a}&=1-v_{fr;f}
\end{aligned}
\end{equation}

\begin{figure}
\centering
\includegraphics[width=0.45\textwidth]{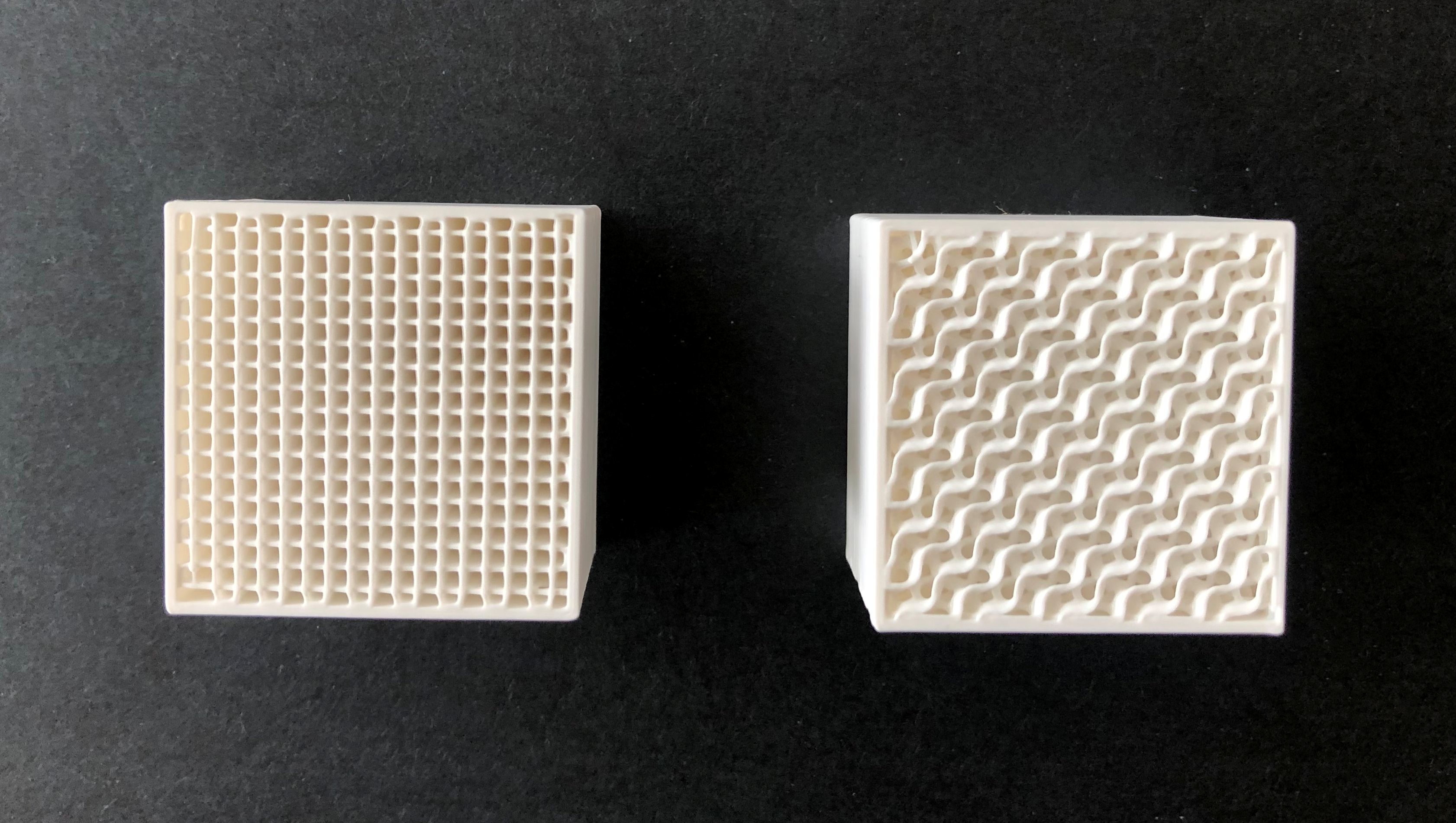}
\caption{Examples of FFF printed specimens with rectilinear (left) and gyroid (right) infill patterns (both at 25\% infill density)}
\label{fig:infill_type}
\end{figure}

\indent The samples for which the infill density and pattern were varied were printed with an extrusion factor of 1.0 and an infill density $\leq$100\%. These samples were printed in the shape of a block with the dimensions of 30x30x20 mm. Two infill patterns were chosen; a rectilinear and a gyroid infill pattern (figure \ref{fig:infill_type}). The total steady-state heating time for these specimens was 40 minutes.\\
\indent An overview of all the FFF printed samples with their dimensions, infill characteristics and heating time is given in table \ref{tab:samples_geom}. The additional parameters required to describe the geometry of the cross-section of the printed specimens with varying air void size are listed in table \ref{tab:geo_prop}. 

\begin{table*} 
\centering
\begin{adjustbox}{max width=\textwidth}
\begin{tabular}{cccccc}
\toprule 
\textbf{Sample} & $\mathbf{l\times{w}\times{h}}$ [mm] & \textbf{Infill pattern} & \textbf{Extrusion factor} [-] & \textbf{Infill density} [\%]  & \textbf{Heating time} [min]\\
\midrule
S1 & 30x30x20 & rectilinear & 1.0 & 100 & 40\\
S2 & 30x30x20 & rectilinear & 1.0 & 50 & 40 \\
S3 & 30x30x20 & rectilinear & 1.0 & 25 & 40 \\
S4 & 30x30x20 & gyroid & 1.0 & 50 & 40 \\
S5 & 30x30x20 & gyroid & 1.0 & 25 & 40 \\
S6 & 30x30x2 & rectilinear & 1.0 & 100 & 5 \\
S7 & 30x30x2 & rectilinear & 1.1 & 100 & 5 \\
\bottomrule 
\end{tabular}
\end{adjustbox}
\caption{Overview of the printed samples and their characteristics}
\label{tab:samples_geom}
\end{table*}

\begin{table*}
\centering
\begin{tabular}{ccccccc}
\toprule
\textbf{EF} & $\mathbf{m}$ [g] & $\mathbf{V_{tot}}$ [cm$^3$] & $\mathbf{V_{extr}}$ [cm$^3$] & $\mathbf{v_{fr;f}}$ [-] & $\mathbf{v_{fr;a}}$ [-] & $\mathbf{a}$ [-]\\
\midrule
1.0 & 2.10 & 1.78 & 1.69 & 0.95 & 0.05 & 0.16\\
1.1 & 2.32 & 1.89 & 1.87 & 0.99 & 0.01 & 0.07\\
\bottomrule
\end{tabular}
\caption{Geometrical paramaters of the cross-section of samples printed with different extrusion factors}
\label{tab:geo_prop}
\end{table*}

\subsection {Thermal analysis}
\label{subsec:thermal}
The heat transfer can be divided into various heat exchange modes \cite{Costa2015}: 

\begin{enumerate}
\item \textit{Convection} with the environment,
\item \textit{Radiation} with the environment and between adjacent filaments, 
\item \textit{Conduction} with the printing bed and between adjacent filaments. 
\end{enumerate}

\noindent Generally, conduction is well described by utilizing conductivity parameters found either in literature or provided by filament manufacturers. Filament cooling due to radiative heat exchange between the filaments is negligible \cite{Costa2015}. Radiation with the environment can have an influence on filament temperature when the value of the heat transfer coefficient, the parameter which expresses convection, is relatively low (5 W/m$^{2}\cdot$K) \cite{Costa2015}. However, in most practical applications this value is much higher, which means that overall filament cooling becomes convection controlled \cite{Costa2015}. Thus, in this work heat transfer by radiation was neglected and the focus was on determining the heat transfer coefficient and the temperature of the ambient air that describes the convection with the environment.\\ 
\indent The temperature field $T(\mathbf{x},t)$ is described by the heat equation:
\begin{align}
\centering
\label{eq:heat_eq}
\rho{c_p}\frac{\partial{T(\mathbf{x},t)}}{\partial{t}}=\nabla\cdot(K_{0}\nabla{T(\mathbf{x},t)})+q
\end{align}
\noindent where $c_{p}$ [J/kg$\cdot$K] is the specific heat capacity, $\rho$ [kg$/$m$^{3}$] is the material density, $K_{0}$ [W/m$\cdot$K] is the conductivity of the material, and $q$ [W/m$^3$] is the internal heat source. The thermal properties of PLA shown in table \ref{tab:mat_prop} were used here. The initial temperature of the specimens equaled the room temperature $T_{a}$, thus the initial condition is expressed as:

\begin{equation}
\label{eq:IC}
T(\mathbf{x},0)=T_{a} \quad  \mathbf{x}\in{\Omega}
\end{equation}

\noindent where $\Omega$ represents the domain of the printed specimen. Distinction is made between the boundary conditions at the interface between the heated printing bed and the sample $\Gamma_b$, and at the free surfaces of the sample $\Gamma_f$. The temperature of the heated printing bed is applied as a Dirichlet boundary condition:

\begin{equation}
\label{eq:BC_z0}
T(\mathbf{x},t)=T_{b}\quad \mathbf{x}\in{\Gamma_b}
\end{equation}

\noindent At the free surfaces, it is assumed that the heat exchange between the sample and the environment is governed by convection. The Neumann boundary conditions are expressed as:

\begin{equation}
\begin{aligned}
\label{eq:BC_conv}
K_{0}\frac{\partial{T(\mathbf{x},t)}}{\partial{\mathbf{n}}}+q_c&=0 \quad \mathbf{x}\in{\Gamma_f}\\
q_c&=h(T(\mathbf{x},t)-T_c)
\end{aligned}
\end{equation}

\noindent where $h$ [W/m$^{2}\cdot$K] is the heat transfer coefficient, and $T_c$ is the temperature of the ambient air. These two parameters were calibrated by fitting the numerical simulations to the experimental data acquired in the thermal measurements. Due to the nature of the thermal measurements, it was not entirely clear how the temperature was distributed in the vicinity of the outer surface of the sample. Since the printed specimens were relatively small, the heated printing bed might have influenced the temperature of the air surrounding the samples. This effect was taken into account during the calibration of the heat transfer coefficient. The process was structured in such a way that three situations were taken into account for the temperature of the ambient air:

\begin{enumerate}
\item Case 1\\
It was assumed that the printing bed did not have a significant influence on the temperature of the ambient air. The temperature of the ambient air equaled room temperature $T_{a}$. 
\item Case 2\\
The temperature of the ambient air was higher than the room temperature due to the heated printing bed. The same elevated temperature was assumed for all of the free surfaces of the heated specimen. 
\item Case 3\\
In this case, it was still assumed that the free surfaces on the sides of the specimen, $\Gamma_{f;s}$ were exposed to air of which the temperature was higher than room temperature. However, since the top surface of the specimen, $\Gamma_{f;t}$ was more distant from the printing bed, a different temperature was assumed for this surface.
\end{enumerate}

\subsection{Numerical set-up}
\label{subsec:num_setup}
\begin{figure*}
\centering
\includegraphics[width=0.75\textwidth]{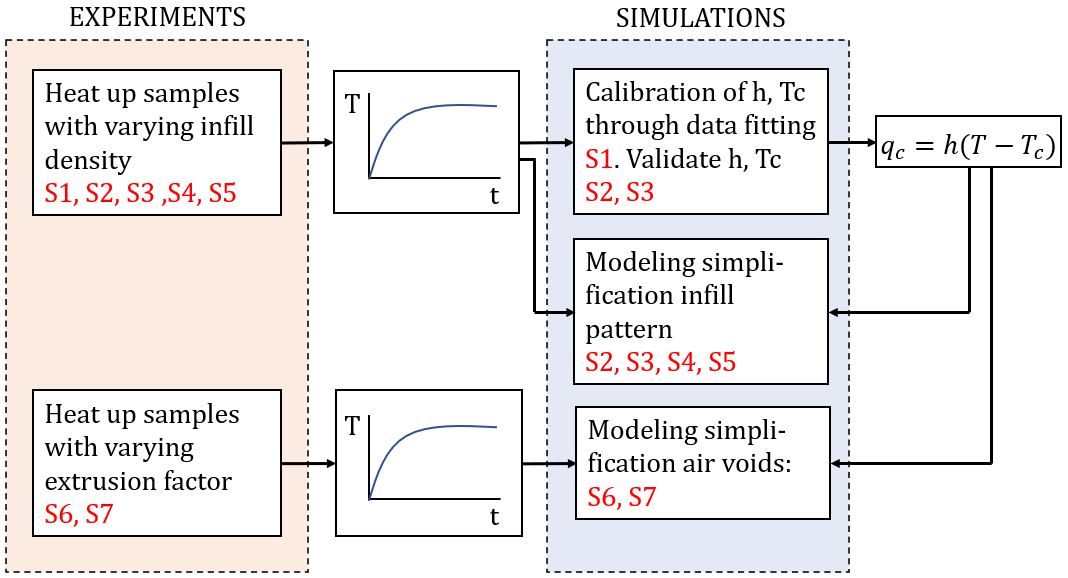}
\caption{Summary of experimental and numerical methodology used in this paper}
\label{fig:method}
\end{figure*}

All finite element (FE) analyses simulating the heat transfer in the printed specimens were set up in Ansys (Mechanical APDL 19.2). The analysis type was a transient thermal analysis. The initial and boundary conditions are listed in equations \ref{eq:IC}-\ref{eq:BC_conv}. The FE samples were discretized with 8-node thermal finite elements (SOLID70). The default element size coincided with the width of the extruded filament and layer height of the printed specimens. This is the case unless it is explicitly stated that a different mesh discretization was used. Convergence analyses were performed prior to all simulations presented in this paper to ensure that the chosen mesh sizes were appropriate. As the air content in the printed samples was included in some simulations, distinction was made between the material properties of the air and of the PLA. Both materials were discretized with the same element type. All analyses that were performed in this work were assumed to be physically linear. This assumption is deemed valid as the maximum temperature that the specimens were heated up to is lower than the glass transition temperature of PLA. Thus, the constant thermal properties as given in table \ref{tab:mat_prop} were used in the numerical simulations. In the simulations, the measured value of 56$^{\circ}$C was used for the temperature of the heated printing bed, unless it is explicitly stated otherwise.\\
\indent The overall experimental and numerical methodology is summarized in figure \ref{fig:method}. On the left side of the chart it is shown that the influence of the infill density (and pattern) is measured experimentally in specimens S1 - S5. Specimens S6 and S7 are used for the measurement of the influence of the air voids. Prior to numerical validation of these two factors, the convective thermal parameters $h$ and $T_c$ are numerically calibrated by using the experimental measurements on specimen S1, and numerically validated by comparing numerical and experimental results on specimens S2 and S3.

\section{Experimental results} 
\label{sec:exp_res}

\subsection{Varying infill density}
\label{subsec:ID}

\begin{figure*}
\centering
\subfloat[100\% infill density]{\includegraphics[width=0.5\textwidth]{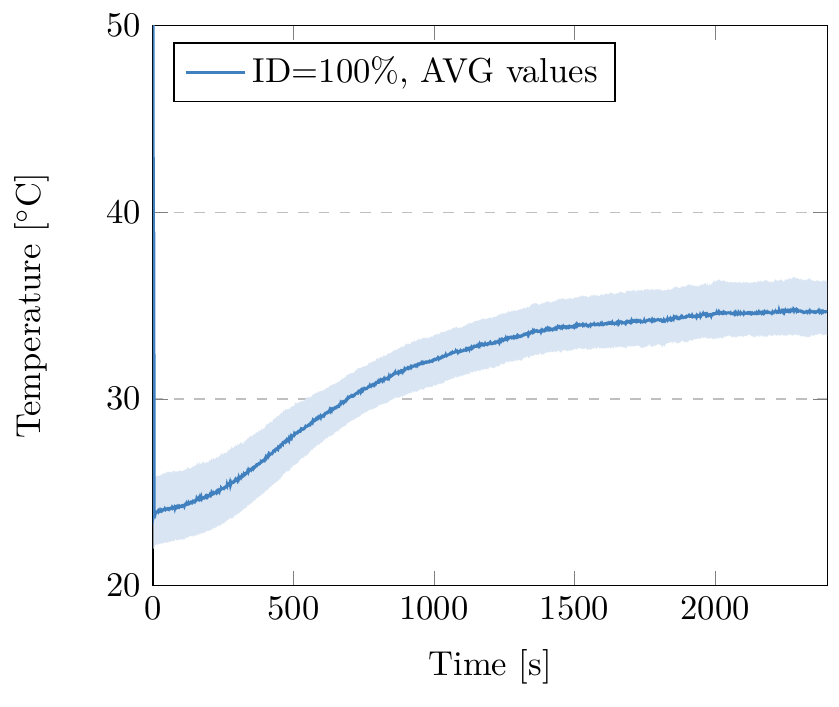}}\hfill
\subfloat[50\% infill density]{\includegraphics[width=0.5\textwidth]{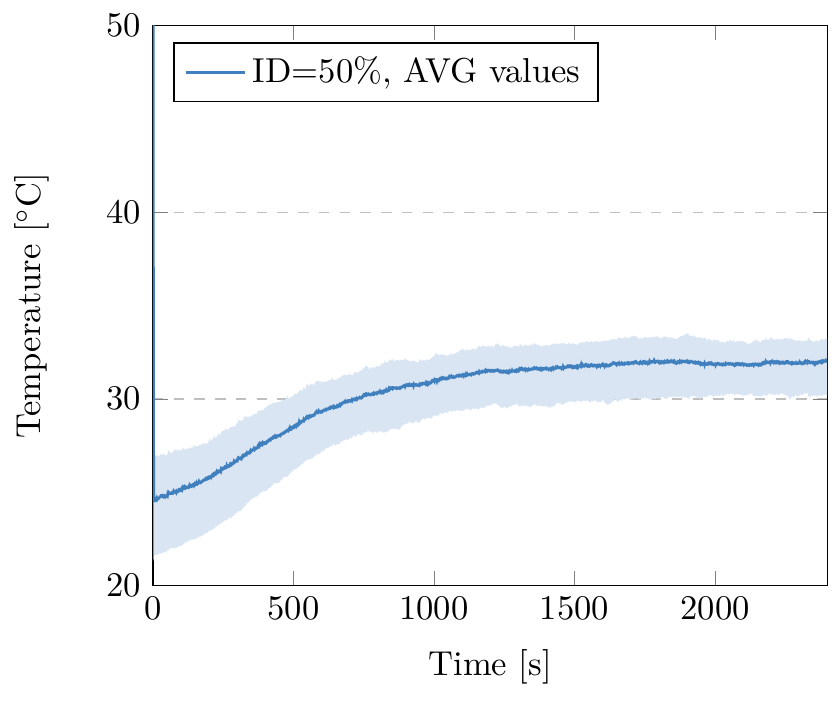}}\\
\subfloat[25\% infill density]{\includegraphics[width=0.5\textwidth]{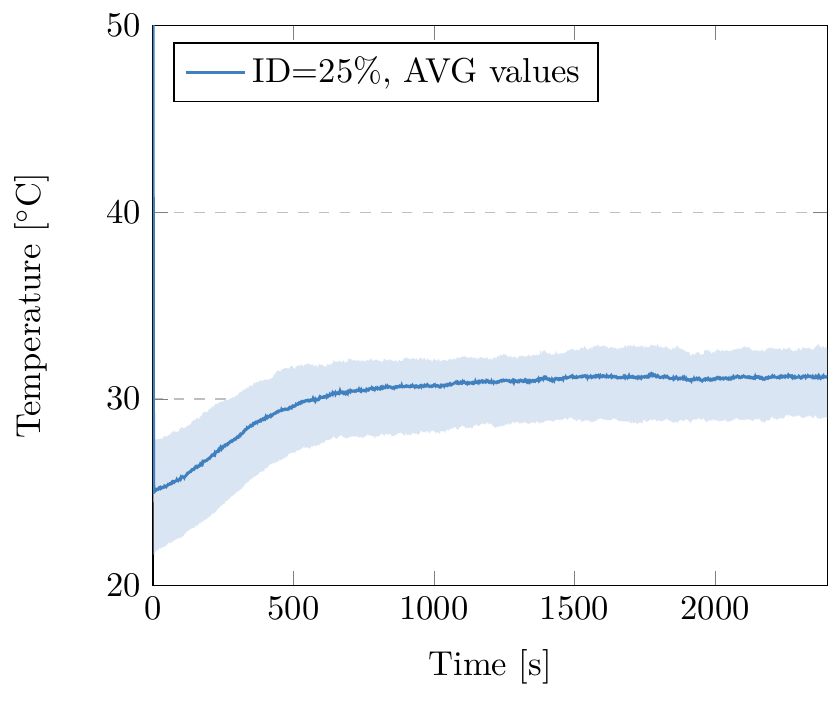}}\hfill
\subfloat[Influence infill density]{\includegraphics[width=0.5\textwidth]{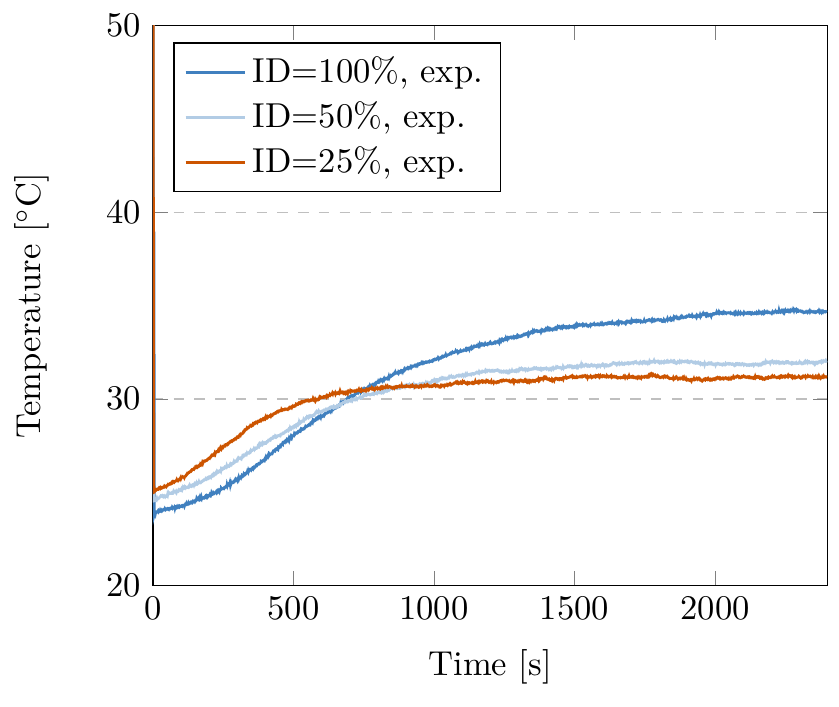}}
\caption{Experimental results of thermal measurements samples with varying infill density}
\label{fig:exp_ID}
\end{figure*}

Specimens S1, S2 and S3 were heated to capture the influence of the infill density on the heat flow. Figure \ref{fig:exp_ID} shows the measured average temperature and envelope for each specimen as well as a comparison of the measured average temperature between the specimens with different infill density. The envelope data consists of three measurement repetitions on five data points for each measurement (15 data points total). It can be seen that the infill density affects the steady-state temperature and the heating rate. The significant increase in air content for the samples with lower infill density results in less mass to be heated, which explains the faster heating rate.

\subsection{Varying extrusion factor}
\label{subsec:EF}

\begin{figure*}
\centering
\subfloat[EF=1.0]{\includegraphics[width=0.5\textwidth]{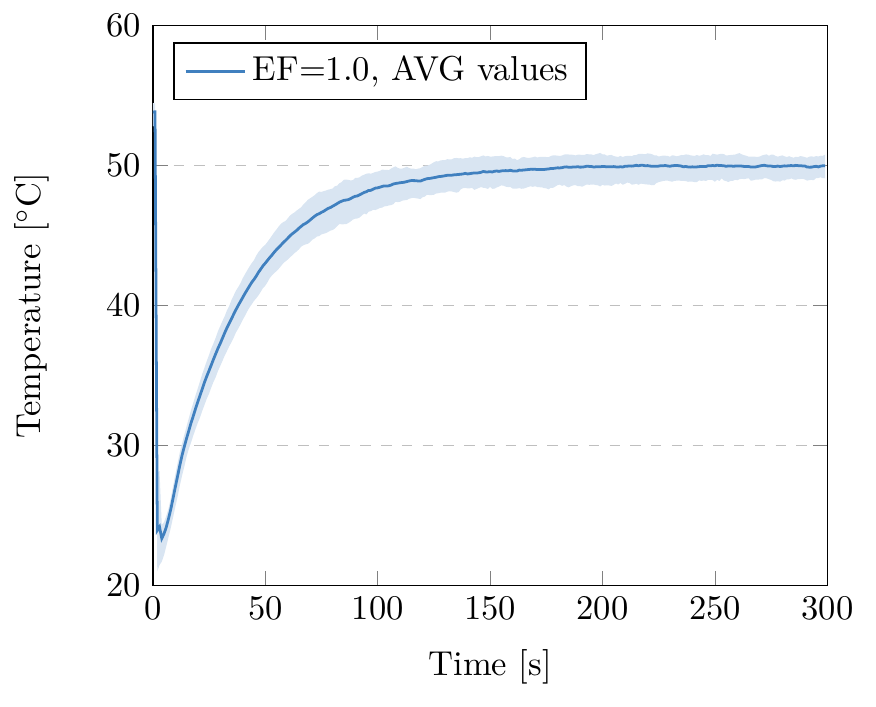}}\hfill
\subfloat[EF=1.1]{\includegraphics[width=0.5\textwidth]{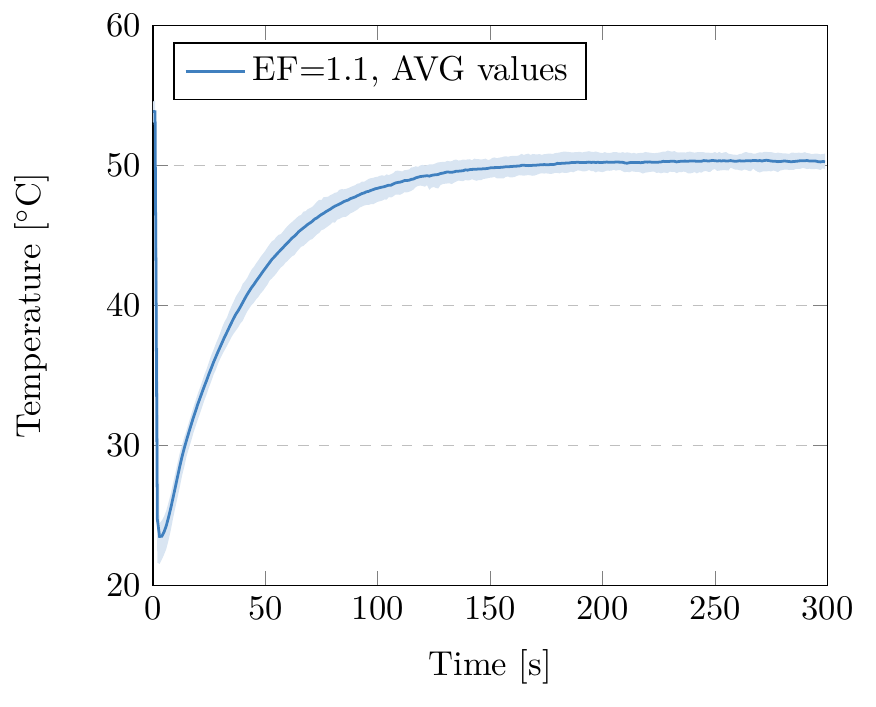}}\\
\subfloat[Influence extrusion factor]{\includegraphics[width=0.5\textwidth]{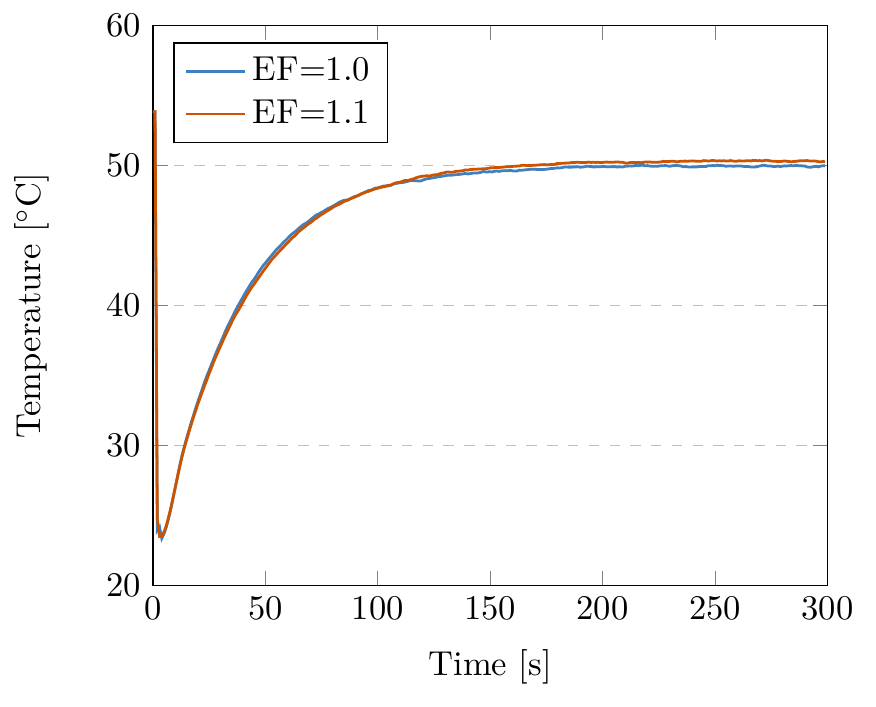}}
\caption{Results thermal measurements on specimens printed with 100\% infill density and varying extrusion factors}
\label{fig:exp_EF}
\end{figure*}

Figure \ref{fig:exp_EF} shows the results of the thermal measurements on the specimens S6 and S7 with varying air void size. Both the average values and the envelopes are plotted. Additionally, a comparison is made between the results of samples S6 and S7 to see what the influence is of the air void content. It can be seen that there is no significant difference between the temperatures measured on the samples printed with different extrusion factors. Varying the extrusion factor and thus the air void size does not seem to influence the heat transfer in the printed samples.   

\subsection{Varying infill pattern}
\label{subsec:IP}

\begin{figure*}
\centering
\subfloat[Infill density=50\%]{\includegraphics[width=0.5\textwidth]{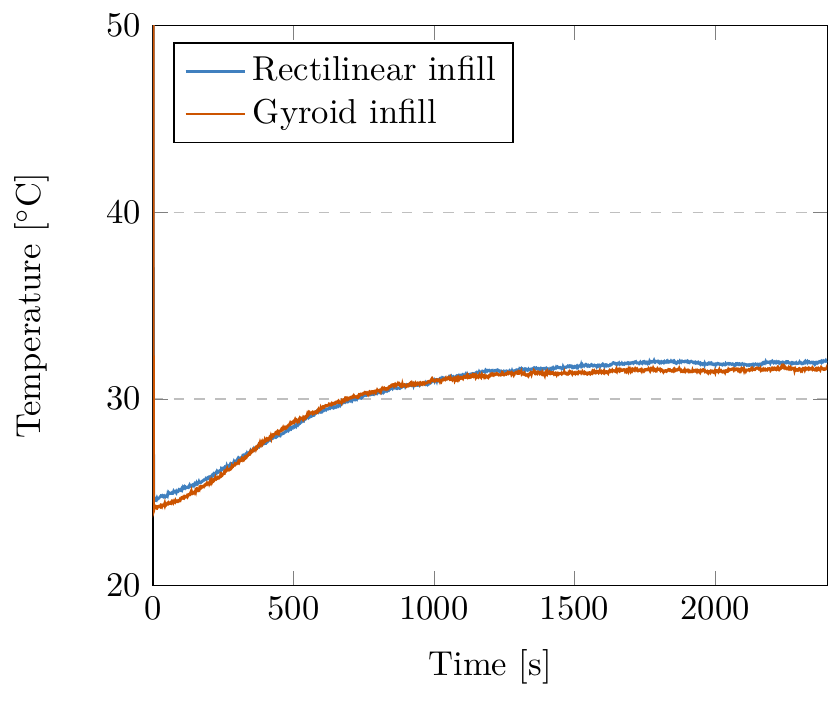}}\hfill
\subfloat[Infill density=25\%]{\includegraphics[width=0.5\textwidth]{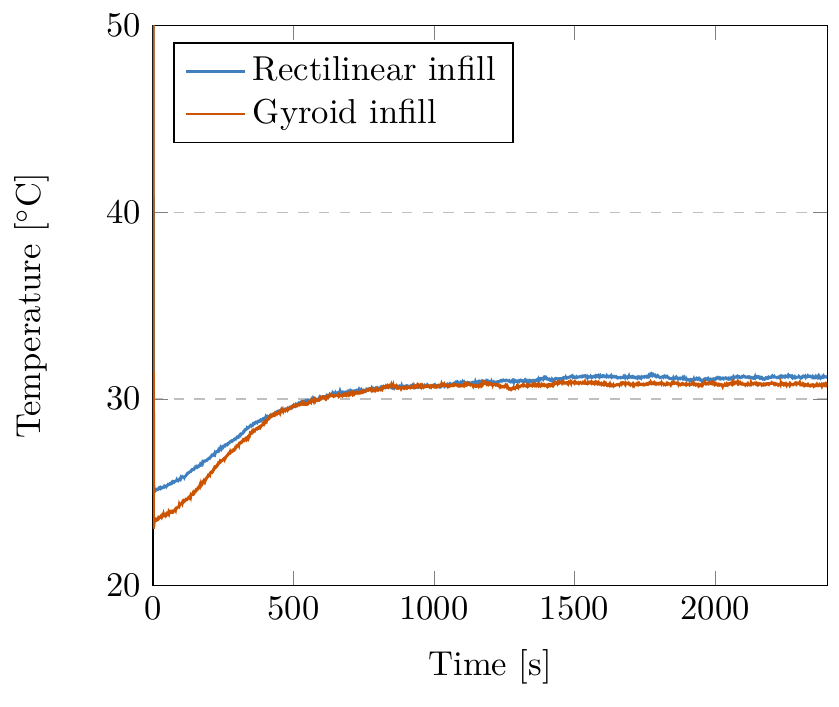}}
\caption{Influence of infill pattern on heat transfer (experimental results)}
\label{fig:RL_GYexp}
\end{figure*}

To see whether the type of infill pattern influences the heat transfer, a comparison is made between samples printed with a rectilinear infill pattern and samples with a gyroid infill pattern. The comparison is done at an infill density of 50\% (S2 and S4) and 25\% (S3 and S5). The experimental results of the thermal measurements are shown in figure \ref{fig:RL_GYexp}. At both infill densities, the differences between the measured temperatures of the rectilinear and gyroid specimen are negligible. Both the heating rate and steady-state temperature are similar for the samples with a rectilinear and gyroid infill pattern. 

\section{Numerical results \& discussion}
\label{sec:num_sim}

\subsection{Determining the convective parameters}
\label{subsec:fitting}

The three cases as described in section \ref{subsec:thermal} were used in the calibration of the convective parameters, $h$ and $T_c$. The numerical results presented below were fitted to the the experimental results from the thermal measurements on sample S1. The experimental thermal measurements of samples S2 and S3 were used for validation of the determined parameters.

\subsubsection{Case 1}
\begin{figure*}
\centering
\includegraphics[width=0.5\textwidth]{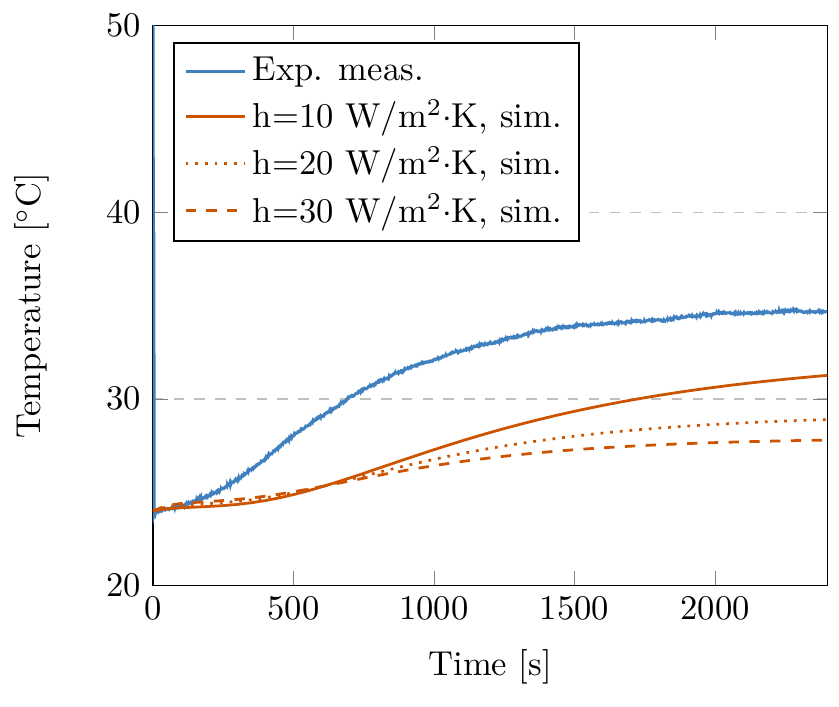}
\caption{Numerical simulations compared to experiments (case 1)}
\label{fig:ID100_h}
\end{figure*}

The first situation considered in the fitting process was the one in which the ambient temperature was assumed to be equal to the room temperature of 25$^{\circ}$C. This value was also used as the initial temperature of the FE sample in the numerical simulations. The heat transfer coefficient in the simulations was varied between 10 and 30 W/m$^{2}\cdot$K. The results are plotted in figure \ref{fig:ID100_h} for $h$=10, 20, 30 W/m$^{2}\cdot$K. It is obvious that none of the simulations show agreement with the measurements. In each of the simulations, the converged temperature is significantly lower than the steady-state temperature found in the measurement (within the time frame of 40 minutes). Increasing the heat transfer coefficient above 30 W/m$^{2}\cdot$K leads to an even lower steady-state temperature and a faster heating rate. For $h$=10 W/m$^{2}\cdot$K, the steady-state is not reached within the considered time span, so using a smaller value for the heat transfer coefficient will delay this even further. 

\subsubsection{Case 2}
\begin{figure*}
\centering
\subfloat[h=20 W/m$^{2}\cdot$K]{\includegraphics[width=0.5\textwidth]{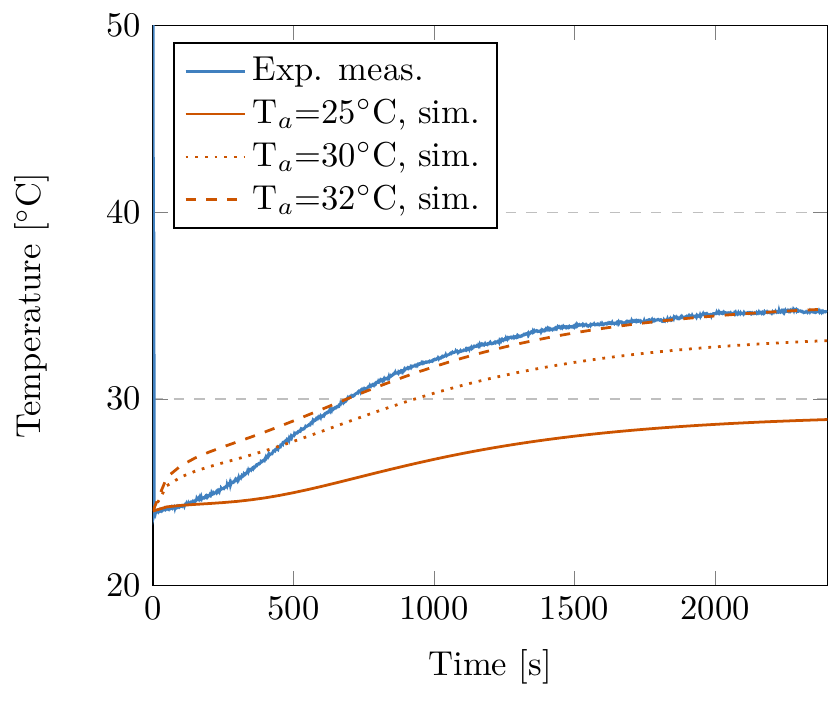}}\hfill
\subfloat[h=30 W/m$^{2}\cdot$K]{\includegraphics[width=0.5\textwidth]{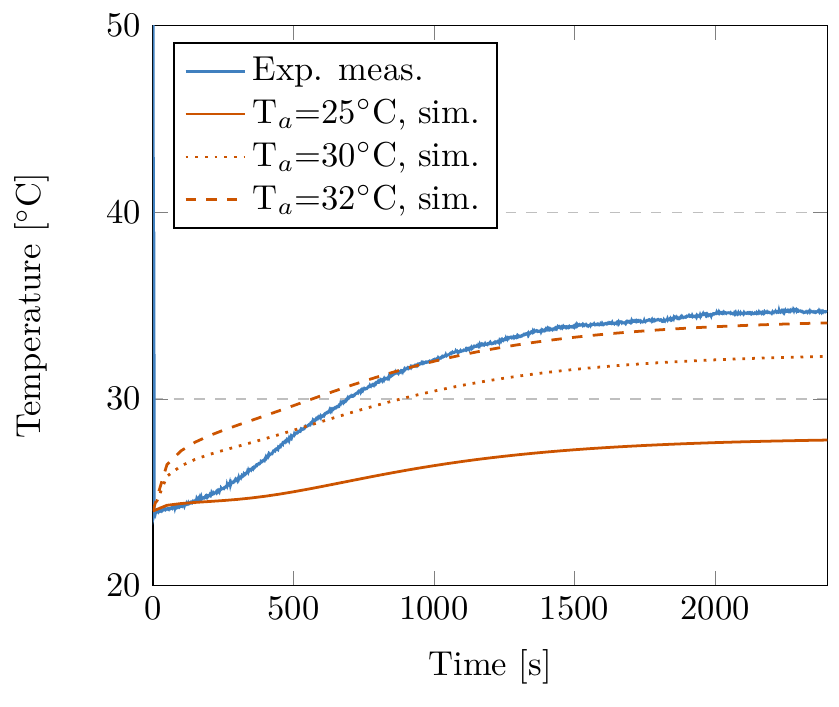}}
\caption{Numerical simulations compared to experiments (case 2)} 
\label{fig:ID100_hTa}
\end{figure*}

The assumption of the printed specimen being subjected to room temperature convection in the vicinity of its free surfaces seems to be incorrect. The simulations were therefore repeated for the same values of $h$ as in figure \ref{fig:ID100_h}, but in combination with higher values for the ambient temperature. In figure \ref{fig:ID100_hTa} the results are plotted for the simulations in which the ambient temperature was varied between 25-32$^{\circ}$C, with heat transfer coeffients of 20 W/m$^{2}\cdot$K and 30 W/m$^{2}\cdot$K. The same convective parameters were applied at all free surfaces. It can be seen that the correct steady-state temperature is reached for $h$=20 W/m$^{2}\cdot$K and $T_c$=32$^{\circ}$C. However, the temperature at the start of the simulation increases much faster than in the experimental measurements. This problem is not solved by choosing a different value for $h$ than those displayed in figure \ref{fig:ID100_hTa}. 

\subsubsection{Case 3}
\begin{figure*}
\centering
\subfloat[$T_{c;t}$=27$^{\circ}$C, $T_{c;s}$=56$^{\circ}$C]{\includegraphics[width=0.5\textwidth]{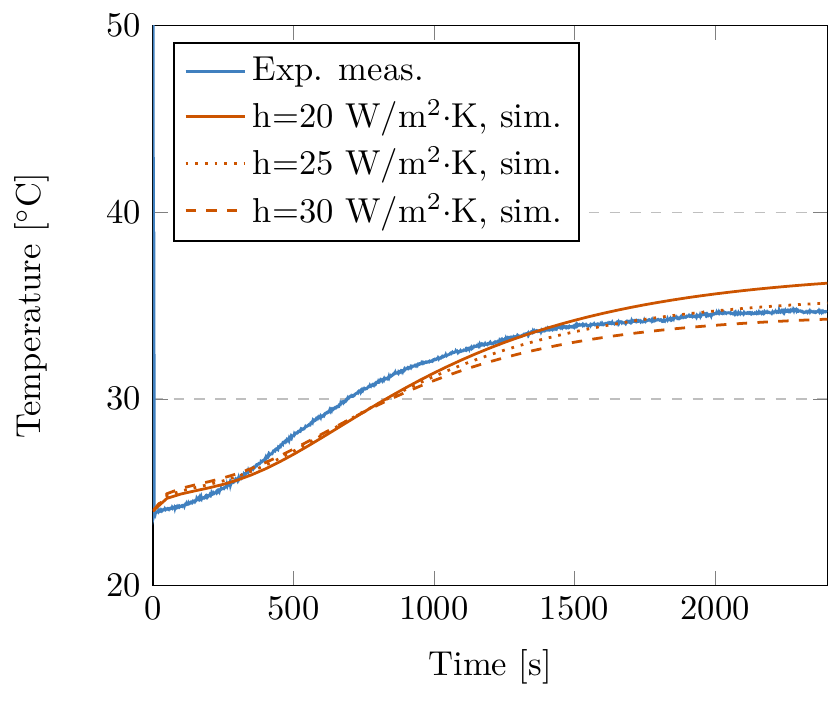}}\hfill
\subfloat[$T_{c;t}$=27$^{\circ}$C, $T_{c;s}$=60$^{\circ}$C]{\includegraphics[width=0.5\textwidth]{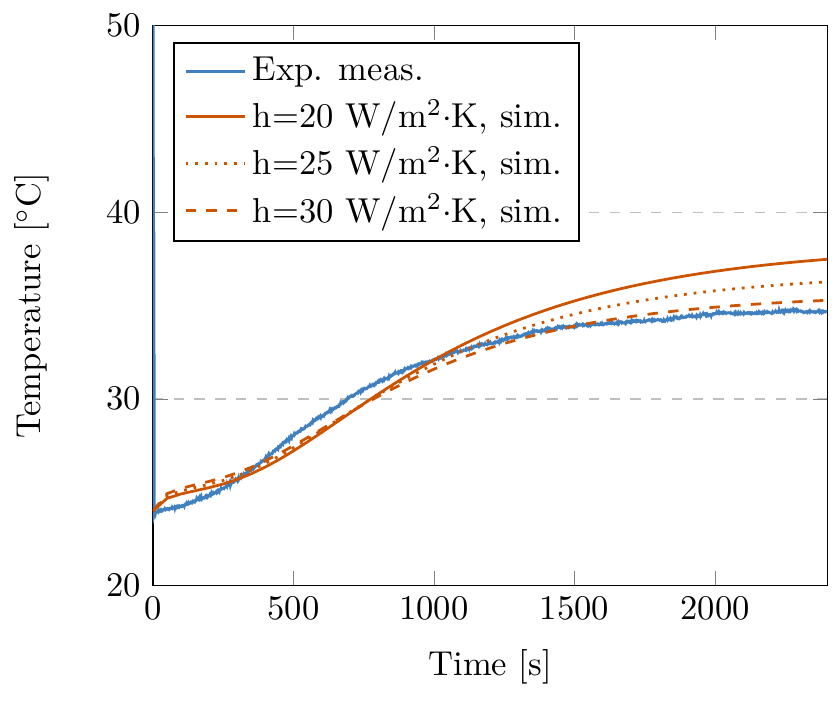}}
\caption{Fitting of numerical simulations to experiments (case 3)}
\label{fig:ID100_hTa2756}
\end{figure*}

In an attempt to tackle the overestimated heating rate at the start of the simulations, the ambient air temperature at the top surface $T_{c;t}$, was varied between 25$^{\circ}$C and 32$^{\circ}$C. The temperature at the side surfaces of the specimens $T_{c;s}$, was chosen to be higher than room temperature, again under the assumption that the printing bed heated up the air above it too. This temperature was varied between 32$^{\circ}$C and the bed temperature. Heat transfer coefficients of 20, 25 and 30 W/m$^{2}\cdot$K were prescribed at all free surfaces. In figure \ref{fig:ID100_hTa2756} the results are plotted for the case in which the temperature at the side surfaces is equal to the bed temperature and the temperature at the top surface equals 27$^{\circ}$C. It can be seen that these values for the air temperature combined with a heat transfer coefficient of 25 W/m$^{2}\cdot$K provide a good match with the experimental measurements. \\
\indent Typically, the exact temperature on the surface of the printing bed is not measured but assumed to be equal to the value prescribed in the printing settings. To account for this, the data fitting process was repeated for $T_{b}=T_{c;s}$=60$^{\circ}$C. These results are shown in figure \ref{fig:ID100_hTa2756}. It can be seen that $h$=30 W/m$^{2}\cdot$K combined with $T_{c;s}$=60$^{\circ}$C and $T_{c;t}$= 27$^{\circ}$C also provides a good fit with the experimental results. In the remainder of this paper, the following convective parameters are used: $h$=25 W/m$^{2}\cdot$K, $T_{c;s}$=56$^{\circ}$C, $T_{c;t}$=27$^{\circ}$C.\\

\subsubsection{Validation $h$ \& $T_c$}
\begin{figure*}
\centering
\subfloat[Infill density=50\%]{\includegraphics[width=0.5\textwidth]{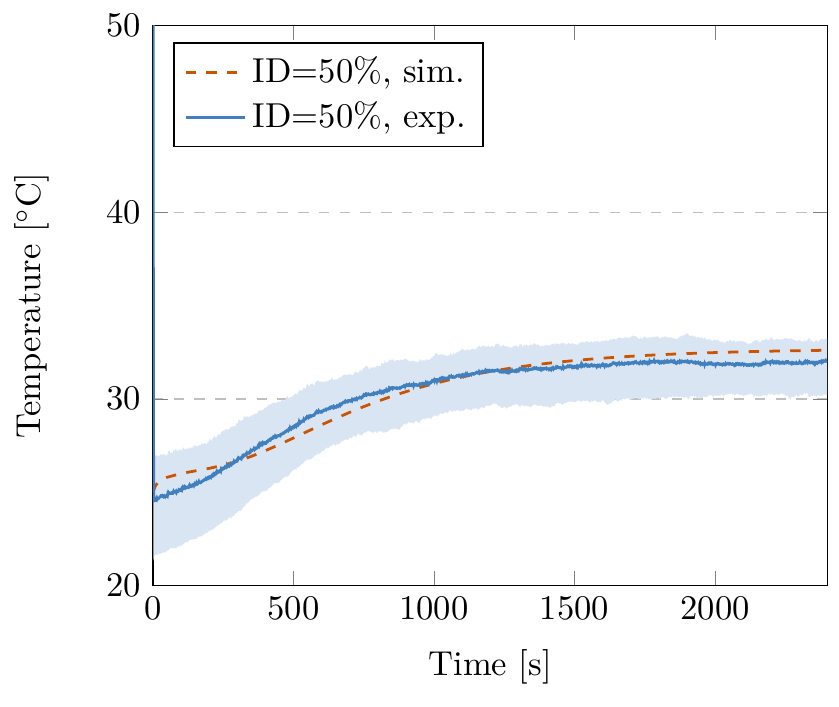}}\hfill
\subfloat[Infill density=25\%]{\includegraphics[width=0.5\textwidth]{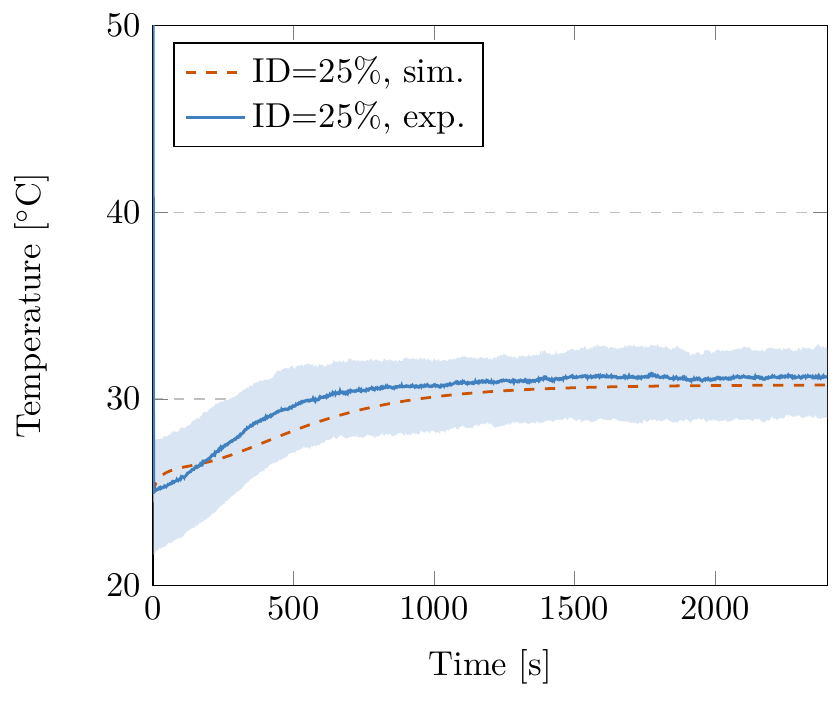}}
\caption{Validation of determined convective boundary conditions: numerical simulations against experiments} 
\label{fig:ID50_ID25}
\end{figure*}

The accuracy of the calibrated convective parameters was validated by simulating the thermal measurements performed on the specimens printed with 50\% and 25\% infill density. The results of the numerical simulations are plotted against the experiments in figure \ref{fig:ID50_ID25} for both of the specimens. Good agreement is found between the simulations and experiments. 

\subsection{Modeling simplifications of the infill structure}
\label{subsec:num_infill}

\subsubsection{Air voids}
\label{subsubsec:num_voids}

\begin{figure*}
\centering
\subfloat[Validation for EF=1.0: exp. vs. sim.]{\includegraphics[width=0.5\textwidth]{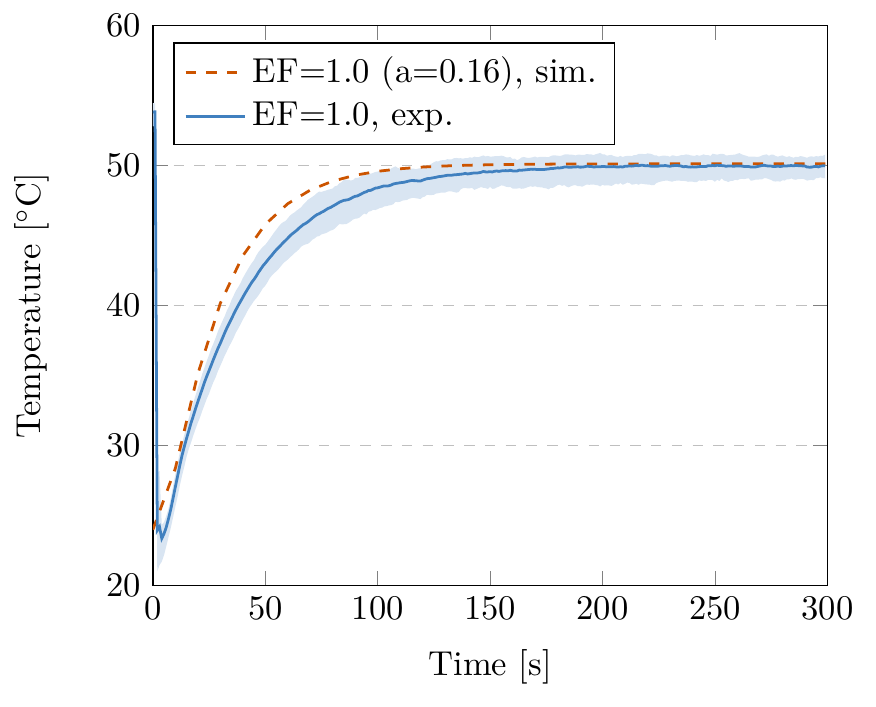}}\hfill
\subfloat[Validation for EF=1.1: exp. vs. sim.]{\includegraphics[width=0.5\textwidth]{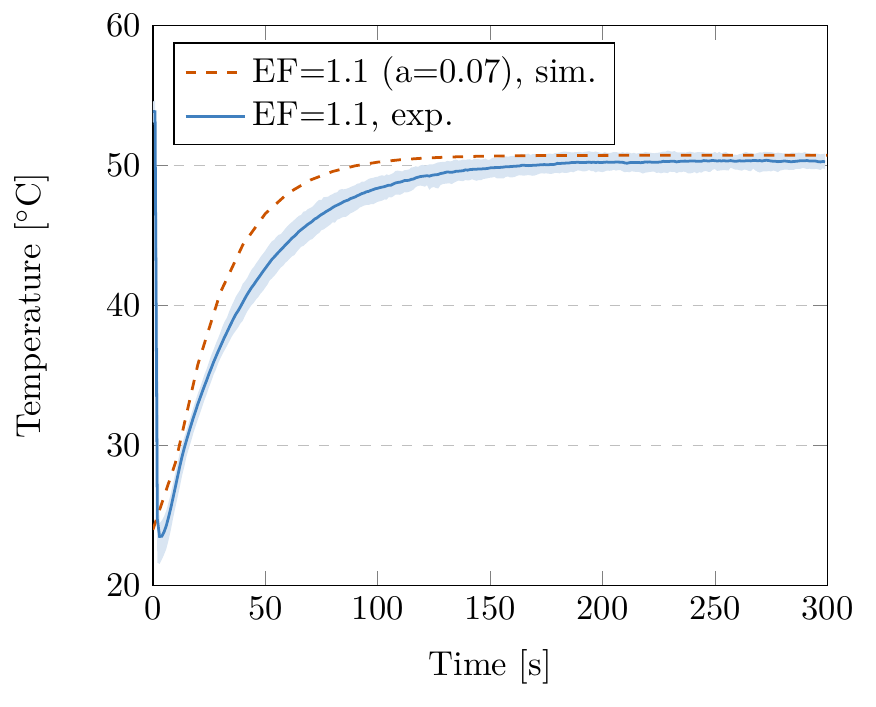}}
\caption{Influence of the air void size: numerical results and experimental validation}
\label{fig:num_voids}
\end{figure*}

\begin{figure*}
\centering
\subfloat[Air voids vs. no air voids (simulations)]{\includegraphics[width=0.5\textwidth]{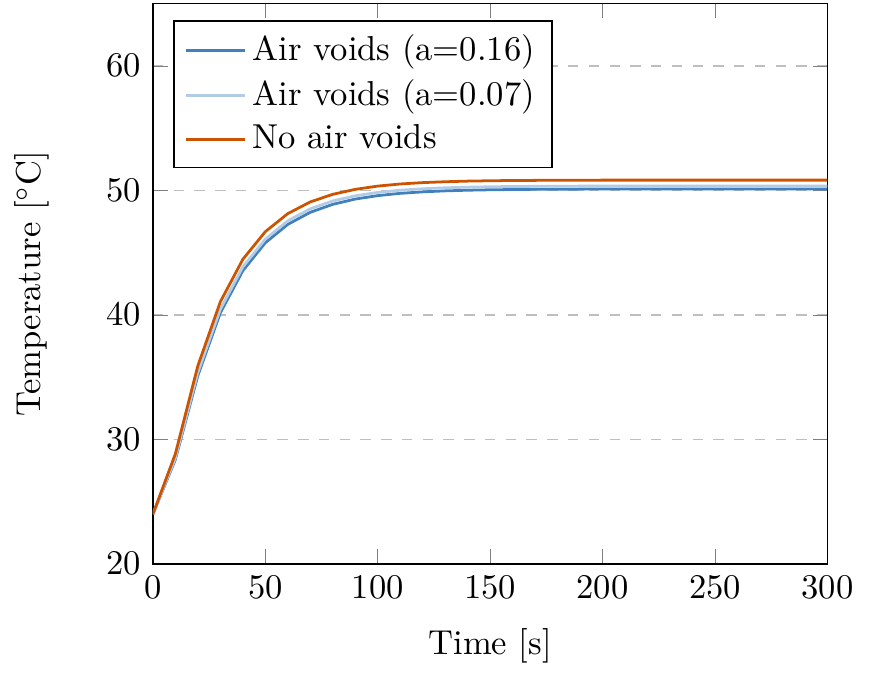}}\hfill
\subfloat[Influence mesh coarseness]{\includegraphics[width=0.5\textwidth]{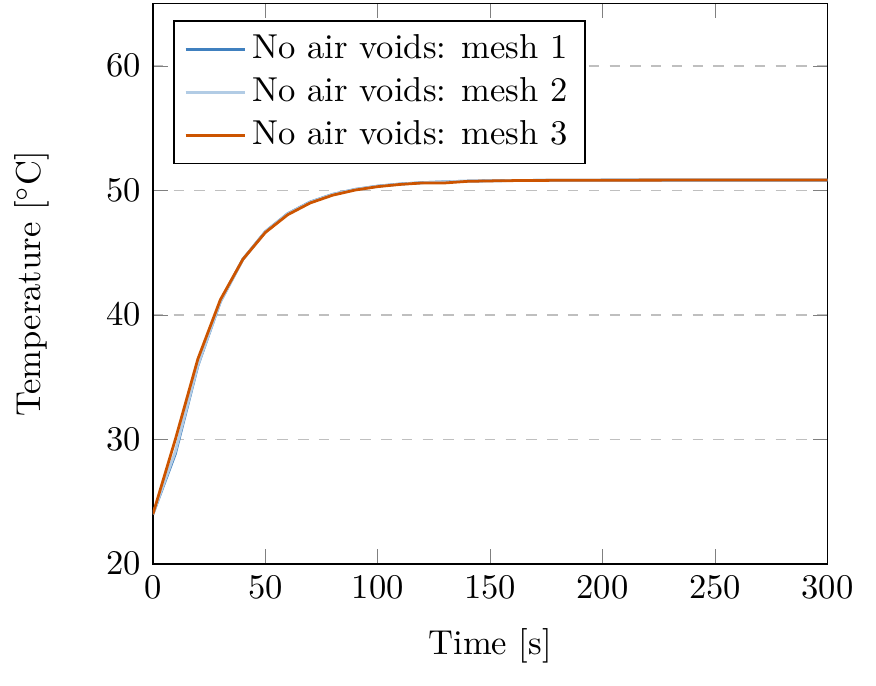}}
\caption{Influence of the air voids and mesh coarseness (simulations)}
\label{fig:num_voids2}
\end{figure*}

In the simulations of the samples with an infill density of 100\% and an extrusion factor of 1.0 presented in section \ref{subsec:fitting}, the FE meshes were made completely dense, ignoring the air voids that were present between the printed filaments. At the same time, the experiments performed on the samples with 100\% infill density and a varying extrusion factor (section \ref{subsec:EF}), implied that role of the air voids was negligible. In this section, we investigate whether a simplified material discretization as used in aforementioned simulations is accurate, or if it is necessary to model the air voids as well. \\
\indent The simulations presented in this section were performed with a FE mesh which didn't include air voids, and with FE meshes which included air voids according to the extrusion factors that were used in the printed specimens (table \ref{tab:geo_prop}). A sample of the detailed mesh including air voids is shown in figure \ref{fig:discretization}. It consists of 13 hexahedral and prismatic elements for the cross-section of a single printed filament with surrounding air compared to one brick element when air voids are not included in the discretization.\\ 
\indent The numerical results are shown in figure \ref{fig:num_voids}. First, the experimental results for the specimens printed with different extrusion factors are compared with the simulations obtained with the meshes including air voids. It can be seen that there is good agreement between the experimental and numerical results. The difference between the steady-state temperatures found in the simulations and experiments is less than 2\%. A comparison is also made between the numerical results for the meshes in which air voids of varying size are included with a mesh in which the material is discretized as a continuum without air voids (figure \ref{fig:num_voids2}). No significant difference is found between the measured temperatures for the different meshes. This coincides with the observation found in the experimental measurements of the air voids size being of insignificant influence on the heat flow. Lastly, the meshes without air voids are coarsened to a mesh where the element dimensions are twice as large as the filament cross-sectional dimensions (mesh 2), and to a mesh where the element dimensions are five times as large as the filament dimensions (mesh 3). As can be seen in figure \ref{fig:num_voids2}, no significant difference is found between the results obtained with these three meshes. 

\subsubsection{Infill density and pattern}
\label{subsubsec:num_IDIP}

\begin{figure*}
\centering
\subfloat[Mesh 1]{\includegraphics[width=0.3\textwidth]{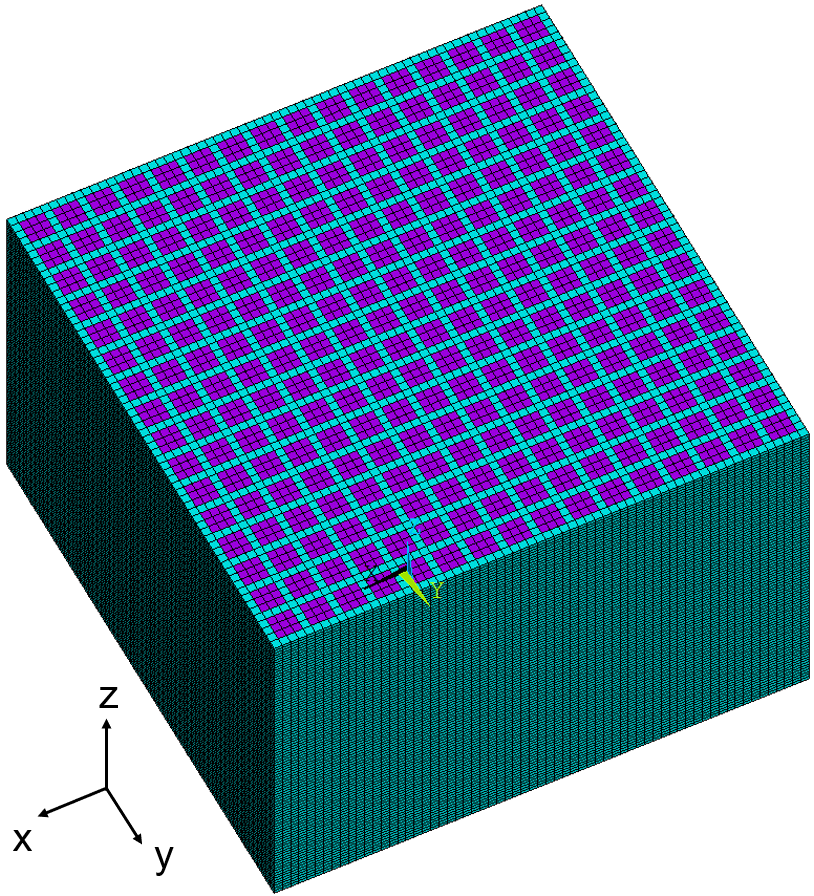}}\hfill
\subfloat[Mesh 2]{\includegraphics[width=0.3\textwidth]{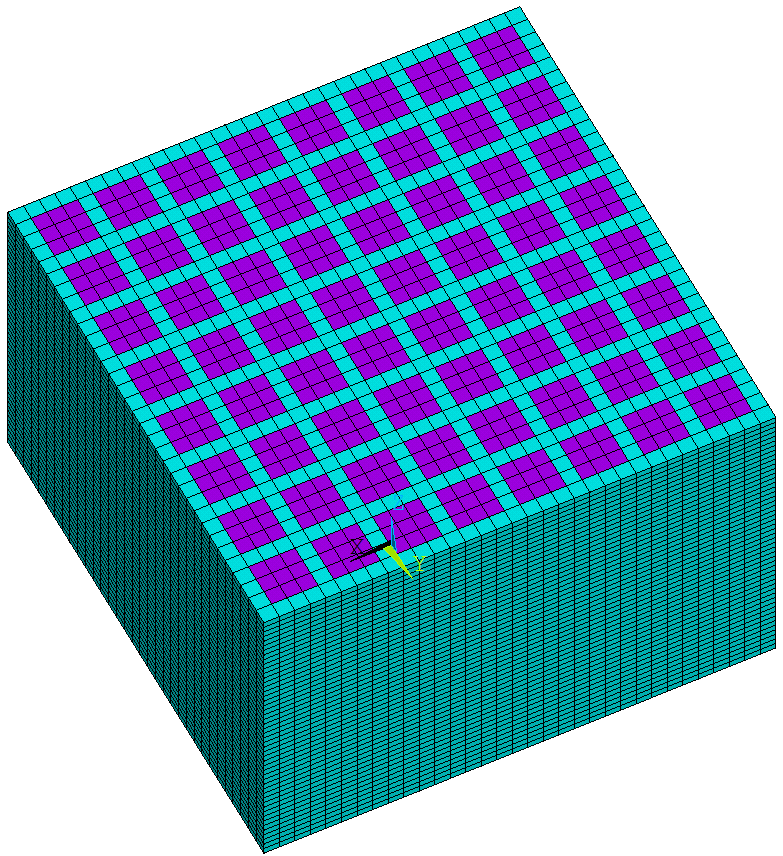}}\hfill
\subfloat[Mesh 3]{\includegraphics[width=0.3\textwidth]{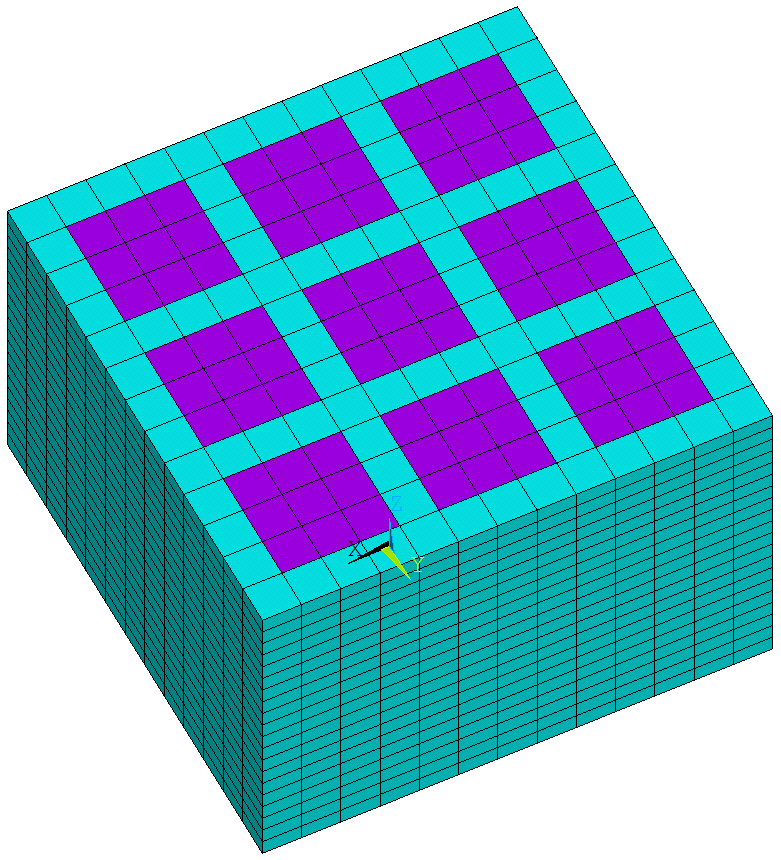}}
\caption{FE meshes of varying coarseness for an infill density of 25\% (blue: PLA elements, purple: air elements)}
\label{fig:ID25_coa}
\end{figure*}

\begin{figure*}
\centering
\subfloat[Numerical results]{\includegraphics[width=0.5\textwidth]{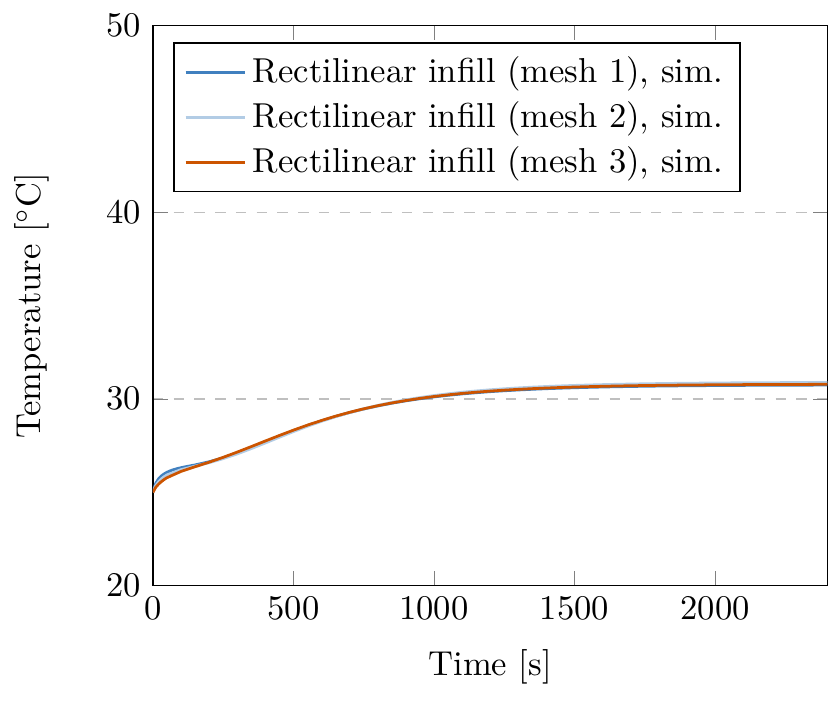}}\hfill
\subfloat[Numerical vs. experimental results]{\includegraphics[width=0.5\textwidth]{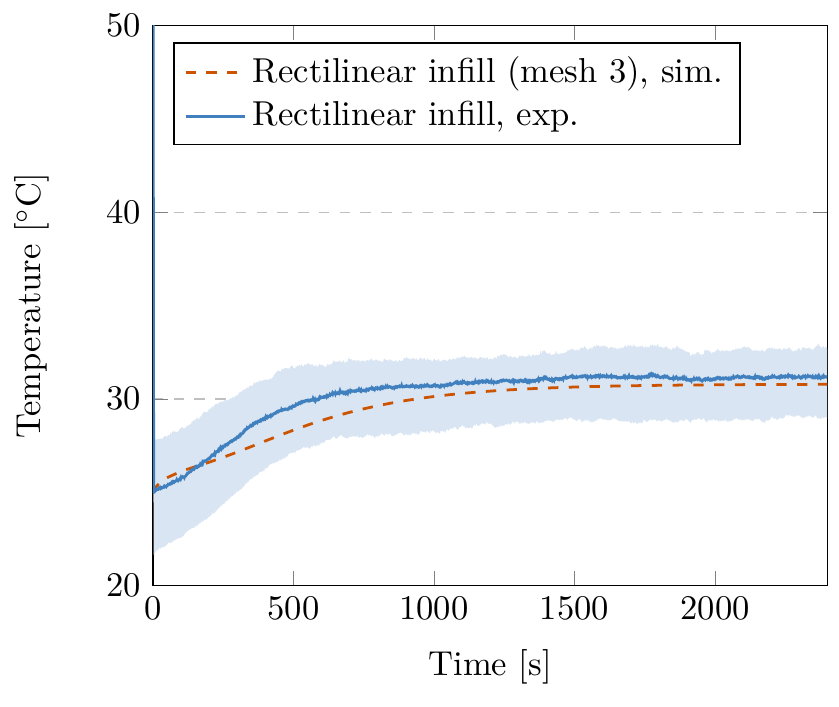}}
\caption{Comparison of various meshes for a rectilinear infill with 25\% infill density}
\label{fig:infl_mesh}
\end{figure*}

\begin{figure*}
\centering
\includegraphics[width=0.5\textwidth]{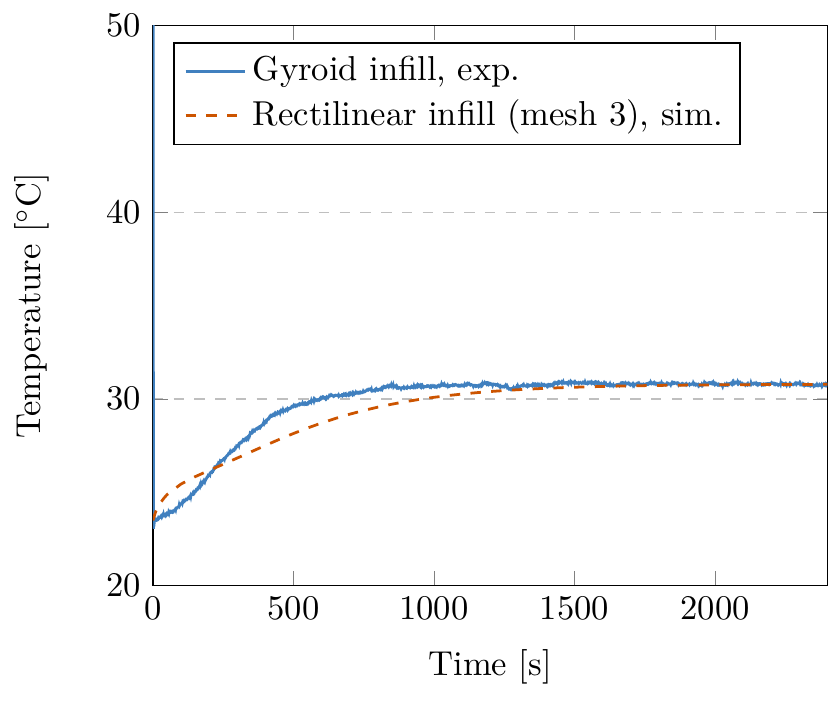}
\caption{Complex infill structure discretized with a simplified coarse mesh: gyroid infill vs. FE mesh 3}
\label{fig:num_exp_GYRL}
\end{figure*}

From section \ref{subsec:ID} we know that varying the infill density has a significant influence on the heating rate and steady-state temperature (figure \ref{fig:exp_ID}), whereas the exact infill pattern doesn't seem to influence the aforementioned characteristics (figure \ref{fig:RL_GYexp}). This means that it might be possible to simulate the accurate heat transfer within complex infill geometries by using much coarser and simplified FE meshes. To test this hypothesis, numerical models were set-up which simulated the heat transfer in the samples printed with a 25\% rectilinear infill. In the models the geometry was discretized with the correct infill density and type, but with meshes of varying coarseness. In the default mesh the element dimensions equaled the dimensions of the filament cross-section. Based on the default mesh, two other, coarser FE meshes were used. In mesh 2, the element dimensions were twice as large as the filament cross-sectional dimensions. In the coarsest mesh (mesh 3) the element dimensions were five times as large as those of the filament cross-section. Since the empty spaces between the printed material were modeled with air elements, the air elements were also coarsened accordingly. The three different meshes are shown in figure \ref{fig:ID25_coa}.\\ 
\indent The results are plotted in figure \ref{fig:infl_mesh}. It can be seen that all of the aforementioned mesh discretizations yield a similar temperature profile. By using the simplified and coarse rectilinear mesh (mesh 3), the computational time was reduced significantly. The number of elements and the CPU time required to solve the numerical models with the various meshes are listed in table \ref{tab:CPU}. When compared to the experimental results of the thermal measurements performed on the printed samples with rectilinear infill, there is also a good agreement.\\
\indent The next modeling simplification was to discretize a complex infill pattern, with both a simplified pattern and a coarser mesh. This has been done by comparing the experimental results of the specimen with a gyroid infill pattern with the coarsest rectilinear finite element model, both with a 25\% infill density. These results are shown in figure \ref{fig:num_exp_GYRL}. Again, there is good agreement between the experiments and the simulations. Both results in figures \ref{fig:infl_mesh} and \ref{fig:num_exp_GYRL} support the experimental findings, namely that exact infill pattern does not influence the temperature profiles, as long as the correct infill density is taken into account.

\begin{table}
\centering
\begin{tabular}{lcc}
\toprule
& Elements [-] & CPU [min]\\
\midrule
Mesh 1 & 422500 & 96\\
Mesh 2 & 54450 & 5\\
Mesh 3 & 3380 & $<$1\\
\bottomrule 
\end{tabular}
\caption{Comparison of CPU for various FE meshes}
\label{tab:CPU}
\end{table}

\section{Conclusions \& Outlook}
\label{sec:conclu}

The overall objective of this work was to provide more clarity on two aspects which affect the accuracy and efficiency of heat transfer simulations of FFF printing. On one hand, a closer look has been taken at the prescription of thermal convective boundary conditions. A value for the heat transfer coefficient and ambient temperature, expressing the natural convection between FFF printed material and its environment have been calibrated through numerical data fitting of experimental thermal measurements.\\
\indent On the other hand, simplifications for modeling air-filled and complex infill structures have been investigated. It has been found that accurate heat transfer can be simulated in such structures, as long as the infill density is respected. Discretization of the exact infill geometry does not lead to a significant increase in accuracy of the heat transfer simulations when compared to the experimental data. In densely packed geometries, the printed material can be modeled as a continuum. Discretization of the air voids between the printed filament is not necessary for an accurate heat transfer prediction. \\
\indent The work performed in this paper was limited to the use of PLA and a Prusa i3 MK3 printer. For practical applications it would be of interest to expand the research to other materials such as ABS and to industrial 3D printers. It must also be noted that the influence of the surface roughness of the printed samples on the convective heat transfer was not taken into account. Lastly, the modeling simplifications proposed in this work have only been applied in thermal simulations. The next step is to model the extrusion process in FFF printing and to investigate the influence of the process parameters on the heat transfer and residual stresses.


\clearpage
\bibliography{lit_paper}
\clearpage
\noindent \textbf{Funding}\\
\indent This work was supported by the European Research Council through the H2020 ERC Consolidator Grant 2019 n. 864482 FDM$^2$. The support of AMOS under project 3: Risk management and maximized operability of ships and ocean structures is also gratefully acknowledged.\\
\\
\noindent \textbf{Conflicts of interest}\\
\indent No conflicts of interest were reported by the authors during the writing of this work.\\
\\
\noindent \textbf{Authors' contributions}\\
\indent \textbf{Nathalie Ramos}: Conceptualization, Methodology, Data acquisition and validation, Writing and visualization (original draft). \textbf{Christoph Mittermeier}: Methodology, Writing (review \& editing). \textbf{Josef Kiendl}: Methodology, Writing (review \& editing), Supervision.

\end{document}